\documentclass[12pt]{amsart}
\usepackage{amsmath,amssymb}
\numberwithin{equation}{section}
\theoremstyle{plain}
  \newtheorem{theorem}{Theorem}[section]

\theoremstyle{definition}

\theoremstyle{remark}
 \newtheorem{remark}[theorem]{Remark}

\newcommand{\xpm}[1]{{[x_{#1},x_{#1}^{-1}]}}
\newcommand{\Ch}[1]{\Tr_{{#1}}q^\der}
\newcommand{\rat}[1]{#1}
\newcommand{\daff}[1]{\widetilde{#1}}
\def\St{{\mathrm St}}
\def\u{u}

\def\fwt{\varpi}
\def\GL{{\mathrm{GL}}}
\def\cas{{C}}
\def\emb{\iota}
\def\charge{d}
\def\unip{{\mathbf e}}
\def\H{{\mathbf H}}
\def\daffH{\daff\H_\kappa}
\def\ratH{\rat{\H}_\kappa}
\def\affH{\H^{\mathrm{aff}}}
\def\ratL{\rat{\mathcal L}_\kappa}
\def\daffL{\daff{\mathcal L}_\kappa}
\def\daffst{\daff{\Delta}_\kappa}
\def\ratst{\rat{\Delta}_{\kappa}}
\def\ratC{\rat{\mathcal C}}
\def\tilratC{\rat{{\mathcal T}}}
\def\daffC{\daff{\mathcal C}}
\def\yy{{\F[\underline{u}]}}

\def\der{{\partial}}

\def\alg{{\mathfrak a}}

\def\sp{{S}}

\def\gm{\g}
\def\hg{{\mathfrak{h}}}
\def\affg{{\wh\g}}

\def\affLl{{L}_{\kappa}}
\def\Pg{X_\m}
\def\Vbox{{V}}
\def\g{{\mathfrak{g}}}
\def\C{{\mathbb F}}

\def\Z{{\mathbb Z}}

\def\F{{\mathbb F}}
\def\W{W}

\def\gl{{\mathfrak{gl}}}
\def\h{E}

\def\al{{\alpha}}
\def\alch{\alpha^\vee}
\def\e{{\epsilon}}
\def\ech{\epsilon^\vee}
\def\lm{{\lambda}}

\def\lsm{{\lambda/\mu}}

\def\l{{\ell}}
\def\m{{m}}

\def\Xreg{{A}}
\def\T{{T}}
\def\Tr{{{\mathrm{Tr}}}}

\def\End{{{\mathrm{End}}}\,}

\def\dim{{{\mathrm{dim}}}\,}

\def\bra{{\langle}}
\def\ket{{\rangle}}

\def\wh{\widehat}
\def\+{\mathop{\oplus}}
\def\*{\mathop{\otimes}}
\begin{document}
\noindent
\title[Double Affine Hecke Algebras]{Double affine Hecke algebras,
conformal coinvariants and Kostka polynomials}
\author{Takeshi Suzuki}
\address{Research Institute for Mathematical Sciences\\
Kyoto University, Kyoto, 606-8502, Japan}
\email{takeshi@kurims.kyoto-u.ac.jp}
\maketitle
\begin{abstract}
We study a class of
representations called ``calibrated representations''
of the degenerate double affine
Hecke algebra and those of
the rational Cherednik algebra
of type ${\mathrm{GL}}_n$.
We give a realization of calibrated irreducible
modules as spaces of coinvariants constructed from integrable
modules over the affine Lie algebra $\widehat{\mathfrak{gl}}_m$.
Moreover, we give a character formula of
these irreducible modules in terms of
a level-restricted Kostka polynomials.
These results were conjectured by Arakawa, Suzuki and Tsuchiya
based on
the conformal field theory.
The proofs using recent
results on the representation theory of the double affine Hecke algebras
will be presented in the forthcoming papers.
\end{abstract}

\bigskip
\begin{center}
  {\bf Introduction}
\end{center}

\medskip
The degenerate double affine Hecke algebras and the rational Cherednik 
algebras are trigonometric and rational degeneration
of Cherednik's double affine Hecke algebras (\cite{Ch;daha}) respectively.
In this note,
we present some concrete results concerning
 a certain class of modules over these two algebras
of type $\GL_n$.
The results 
are parallel and we explain those for the rational Cherednik
algebra below.
Let $\F$ be an algebraically closed field of characteristic $0$.
Let $\h$ denote an $n$-dimensional vector
space with the basis $\{\ech_i\}_{i\in[1,n]}$,
where $[1,n]=\{1,2,\dots,n\}$.
Let $\h^*=\oplus_{i=1}^n\F\e_i$
be the dual space of $\h$,
where $\e_i$ denotes the dual vector of $\ech_i$.
The natural pairing is denoted by
$\bra\, |\,\ket : \h^*\times \h\to\F$.
Put $\al_{ij}=\e_i-\e_j$.
Then
$R=\left\{\al_{ij}|\, i,j\in [1,n],\ i\neq j\right\}$
and
$R^+=\left\{\al_{ij}\in R\mid i<j\right\}$
give a set of roots
and a set of positive roots
of type $A_{n-1}$ respectively.
\par
%
Let $\W$ denote the Weyl group associated with the
root system $R$, which acts on $\h$ and $\h^*$ naturally.
Denote by $s_\al$
the reflection in $\W$ corresponding to
$\al\in R$.
We write  $s_{ij}=s_{\al_{ij}}$.
Let $S(\h)$ (resp. $S(\h^*)$)
denote the symmetric algebra of $\h$ (resp. $\h^*$).

Fix $\kappa\in \F\setminus\{0\}$.
Define
the rational Cherednik algebra $\ratH$ of type $\GL_n$
as the associative $\F$-algebra
which is generated by the algebra $S(\h^*),\ \F W$ and $S(\h)$
and subjects to the following
defining relations:
\begin{align*}
&w h w^{-1} =w (h) \ \
(w\in\W,\ h\in \h),\ \ \
w \zeta w^{-1}= w (\zeta)\ \ 
(w\in\W,\ \zeta\in \h^*),\\
&\left[ h, \zeta\right] =\kappa\,
\bra \zeta\mid h\ket
+\sum\nolimits_{\al\in R^+}\bra \al\mid h\ket
\bra \zeta\mid\alch\ket
s_{\al}\ \ \left(h\in\h,\ \zeta\in\h^*\right). 
\end{align*}
As shown in \cite{GGOR}, the isomorphism classes of
the ``highest weight'' irreducible
$\ratH$-modules are parameterized by
the set of partitions of $n$.
Let  $\ratL(\lm)$ denote the
irreducible $\ratH$-module corresponding to the
partition $\lm$.
Define the elements
$\u_i\ (i\in[1,n])$ of $\ratH$
by
$\u_i=\e_i\,\ech_i+\sum\nolimits_{1\leq j<i}s_{ji}$,
which are pairwise commutative.
We focus on highest weight
$\ratH$-modules
 which are
semisimple over $\F[\u_1,\dots,\u_n]$.
Such modules are  called calibrated modules.
Now, let $\kappa\in\Z_{\geq 1}$.
For $\m\in\Z_{\geq1}$,
put $\Lambda^+_{\kappa}(\m,n)=
\{ \lm=(\lm_1,\dots,\lm_\m)\in\Z^\m\mid
\lm_1\geq\lm_2\geq \dots\geq\lm_\m\geq 1,\
 \kappa-\m-\lm_1+\lm_\m\in\Z_{\geq0},\ 
\sum_{i=1}^\m\lm_i=n\}.$ 
It is shown in \cite{S;ratandtrig}
that the assignment $\lm\mapsto \ratL(\lm)$
induces a bijection between the set $\sqcup_{\m\in[1,\kappa]}
\Lambda^+_{\kappa}(\m,n)$ and the
set of isomorphism classes of
the calibrated irreducible $\ratH$-modules.
%

Let $\m\in\Z_{\geq1}$ and put $\gm=\gl_\m(\F)$.
Let $\affg=\gm\*\C[t,t^{-1}]\+\C c$ be
the affine Lie algebra associated with $\gm$,
where $c$ is a central element.
Put $\gm[t]=\gm\otimes\F[t]$, which is a subalgebra of $\affg$.
Denote the Cartan subalgebra of $\gm$ (resp. $\affg$)
 by $\hg$ (resp. $\wh\hg$),
and its dual space by $\hg^*$ (resp. $\wh\hg^*$).
We identify $\hg^*$ with $\F^\m$ via the natural basis of $\hg$.
For
 $\lm\in\hg^*\cong\F^\m$,
let $\affLl^{\dagger}(\lm)$ denote the
lowest weight irreducible module over $\affg$ with
lowest weight $-\lm-(\kappa-\m) c^*\in\wh\hg^*$,
where $c^*$ denotes the dual vector of $c$.
Let $\Vbox=\C^{\m}$ denote the basic representation of $\gm$,
and regard $\Vbox[x]=\Vbox\*\F [x]$
 as a $\gm[t]$-module naturally.
For a lowest weight $\affg$-module $M$, 
we set
\begin{align*}
\ratC(M)&=
(\Vbox[x_1] \otimes\dots\otimes \Vbox[x_n]\* M)/
\g[t] (\Vbox[x_1] \otimes\dots\otimes \Vbox[x_n]\* M).
\end{align*}
We will introduce an action of $\ratH$
on the space $\ratC(M)$
in Section 4, 
and moreover we have

\medskip\noindent
{\bf Theorem} (Theorem~\ref{th;image_rat})
Let $\kappa, \m\in \Z_{\geq1}$
and  $\lm\in \Lambda_{\kappa}^+(\m,n)$.
Then
$\ratC(\affLl^\dagger(\lm))\cong \ratL(\lm)$
as an $\ratH$-module.
\medskip

Put $\der
=\kappa^{-1}(\sum_{i\in[1,n]}\e_i\ech_i+\sum_{i < j}s_{ij})
=\kappa^{-1}\sum_{i\in[1,n]}\u_i\in\ratH$, and
put $\unip=\frac{1}{n!}\sum_{w\in\W}w$.
Then $[\der,\unip]=0$ and
the space $\unip \ratL(\lm)$ admits a spectral decomposition
with respect to the operator $\der$. We have

\medskip
\noindent
{\bf Theorem} (Theorem~\ref{th;char_rat})
Let $\kappa,\m\in\Z_{\geq1}$
and $\lm\in \Lambda^+_{\kappa}(\m,n)$.
Then
  \begin{equation}
    \Ch{\unip\ratL(\lm)}=
\frac{q^{\cas_\lm+\frac{1}{2}n(n-1)}}{\mathop{\Pi}\limits_{i=1}^n
(1-q^i)}
{K}_{\lm\; (1^n)}^{(\kappa-\m)}(q^{-1}).
  \end{equation}
Here $\cas_\lm=\frac{1}{2\kappa}
\sum_{i=1}^\m \lm_i(\lm_i-2i+1)$,
and
${K}_{\lm\; (1^n)}^{(\kappa-\m)}(q)$
stands for the $(\kappa-\m)$-restricted Kostka polynomial
associated with the partitions $\lm$
and $(1^n)=(1,1,\dots,1)$ (see e.g. \cite{SS}).
\medskip

The statements above are first
 conjectured in \cite{AST} for the degenerate double affine Hecke
algebra $\daffH$.
These 
conjectures are proved
using the combinatorial description of
calibrated representations of $\daffH$ given in \cite{SV}.
It is shown in \cite{S;ratandtrig} that
$\ratH$ is a subalgebra of $\daffH$ and the induction functor
$\daffH\otimes_{\ratH}(-)$ gives a nice correspondence
between the category of $\ratH$-modules and that of $\daffH$-modules,
which implies the theorems for $\ratH$.
The details will be presented in forthcoming papers. 
%
%
\section{The affine Lie algebra}
%
%
Throughout this paper, we use the notation
$[i,j]=\{i,i+1,\dots,j\}$ for $i,j\in\Z$.

Let $\F$ be a field of characteristic $0$.
Let 
 $m\in\Z_{\geq 1}$.
Let $\gm$ denote the Lie algebra $\gl_\m(\F)$
consisting of all the $\m\times \m$ matrices over $\F$.
Let $\affg=\gm\*\C[t,t^{-1}]\+\C c$ be
the associated affine Lie algebra
with the central element $c$ and the commutation relation
$
  [a\* t^p,b\* t^q]=[a,b]\* t^{p+q}+\Tr(ab)
p\delta_{p+q,0}c$
for $a,b\in\gm,\ p,q\in\Z$.
Put $\gm[t]=\gm\otimes\F[t]$ and $\gm[t,t^{-1}]=\gm\otimes\F[t,t^{-1}]$,
which are Lie subalgebras of $\affg$.
Let $e_{ij}$ $(i,j\in[1,\m])$ denote the matrix unit in $\gm$ with
 only non-zero entries $1$
at $(i,j)$-th component.
Put $\hg=\oplus_{i=1}^{\m}\F e_{ii}$ and $\wh\hg=\hg \+\C c$,
which give a Cartan subalgebra of $\gm$ and $\affg$ respectively.
Denote the dual space of $\hg$ (resp. $\wh\hg$)
 by $\hg^*=\oplus_{i=1}^{\m}\F e_{ii}^*$
(resp. $\wh\hg^*=\hg^*\oplus\F c^*$),
where $e_{ii}^*$ and $c^*$ are the dual vectors of
$e_{ii}$ and $c$ respectively.
A $\affg$-module
$M$ is said to be of level $\l\in\F$ if
the center $c$ acts as a scalar $\l$ on $M$.
For $\kappa\in\F$ and $\lm\in\hg^*$, let $\affLl(\lm)$
(resp. $\affLl^{\dagger}(\lm)$)  denote the
highest weight (resp. lowest weight) irreducible
$\affg$-module with highest weight $\lm+(\kappa-\m)c^*$
(resp. lowest weight $-\lm-(\kappa-\m) c^*$).
We identify $\hg^*$ with $\F^\m$
via the correspondence $\sum_{i=1}^{\m}\lm_i e_{ii}^*
\mapsto (\lm_1,\dots,\lm_\m)$,
and introduce its subspaces
$\Pg=\Z^\m$,
$\Pg^+=\{(\lm_1,\dots,\lm_\m)\in \Z^\m\mid
\lm_1\geq\lm_2\geq\dots\geq\lm_\m\}$
and
$\Pg^+(\kappa)=\{
(\lm_1,\dots,\lm_\m)\in \Pg^+
\mid \kappa-\m-\lm_1+\lm_m\in\Z_{\geq0}\}.$
Note that $\affLl(\lm)$ and $\affLl^{\dagger}(\lm)$ are integrable for
$\lm\in\Pg^+(\kappa)$,
and that $\Pg^+(\kappa)$ is empty unless $\kappa-\m\in\Z_{\geq 0}$.

Let $\Vbox=\C^{\m}$ denote the basic representation 
of $\gm$.
Put $\Vbox\xpm{}=\Vbox\*\C\xpm{}$ (resp. $\Vbox[x]=\Vbox\*\F [x]$),
and regard $\Vbox\xpm{}$ (resp. $\Vbox[x]$) as
a $\gm[t,t^{-1}]$-module (resp. $\gm[t]$-module)
through the correspondence
$ a\*t^p\mapsto a\* x^p$ 
$(a\in\gm,\ p\in\Z).$ 
\section{The rational Cherednik algebra and 
the degenerate double affine Hecke algebra}
Let $n\in\Z_{\geq 2}.$
Let $\h$ denote the $n$-dimensional vector
space with the basis $\{\ech_i\}_{i\in[1,n]}$.
Introduce the non-degenerate symmetric bilinear form $(\ |\  )$
on $\h$ by
$(\ech_i|\ech_j)=\delta_{ij}$.
Let $\h^*=\oplus_{i=1}^{n}\F\e_i$
be the dual space of $\h$,
where $\e_i$ is  the dual vector of $\ech_i$.
The natural pairing is denoted by
$\bra\, |\,\ket : \h^*\times \h\to\F$.
Put $\al_{ij}=\e_i-\e_j,\, \alch_{ij}=\ech_i-\ech_j$. 
Then
$R=\left\{\al_{ij}|\,i,j\in [1,n],\ i\neq j\right\}$ and
$R^+=\left\{\al_{ij}\in R\mid i<j\right\}$
give a set of roots and a set of
positive roots of type $A_{n-1}$ respectively.
%
Let $\W$
denote the Weyl group
associated with the
root system $R$, which is isomorphic to the symmetric group
of degree $n$.
Denote by $s_\al$
the reflection in $\W$ corresponding to
$\al\in R$.
We write 
 $s_{ij}=s_{\al_{ij}}$ and $s_i=s_{i\, i+1}$.
Put
$P=\+_{i\in[1,n]} \Z\e_i$.
We denote the symmetric algebra of $\h$ (resp. $\h^*$)
by $S(\h)$ (resp. $S(\h^*)$), on which $\W$ acts.

Fix $\kappa\in \F\setminus\{0\}$.
{\it The rational Cherednik algebra} $\ratH$ of type $\GL_n$
is the associative $\F$-algebra
which is generated by the algebras $S(\h^*),\ \F W$ and $S(\h)$,
and subjects to the following
defining relations:
\begin{align*}
&s_i h s_i =s_i (h) \
(i\in[1,n-1],\ h\in \h),\ \
s_i \zeta s_i= s_i (\zeta)\
(i\in[1,n-1],\ \zeta\in \h^*),\\ 
&\left[ h, \zeta\right] =\kappa\,
\bra \zeta\mid h\ket
+\sum\nolimits_{\al\in R^+}\bra \al\mid h\ket
\bra \zeta\mid \alch\ket
s_{\al}\ \left(h\in\h,\ \zeta\in\h^*\right).
\end{align*}
It is known that 
the natural multiplication map gives an isomorphism
$S(\h^*)\* \F W \* S(\h)\cong\ratH$ as a vector space
(\cite{EG}).

For $\lm\in\Z^\m=\Pg$, we write 
$\lm\models n$ when
 $\sum_{i=1}^{\m}\lm_i=n$
and $\lm_i\in \Z_{\geq 0}$ for all $i\in[1,m]$.
Let $\lm\in \Pg^+$ such that $\lm\models n$.
We sometimes identify $\lm$ with the Young diagram
\begin{equation}\label{eq;youngdiagram}
\lm=\{(a,b)\in\Z^2~\mid~ a\in[1,m],\ b\in[1,\lm_a]\}\subset\Z^2.
\end{equation}
Let $\sp_\lm$ denote the irreducible module of $\W$
associated with the partition $\lm$.
Set $\ratst(\lm)=\ratH\otimes_{\F W\cdot S(\h)} \sp_{\lm}$,
where we let $S(\h)$ act on $\sp_{\lm}$ through
the augmentation map given by
$\ech_i\mapsto0\ (i\in[1,n])$.
It is known that the $\ratH$-module $\ratst(\lm)$ has
a unique simple quotient~\cite{GGOR},
which we will denote by $\ratL(\lm)$.

{\it The degenerate double affine Hecke algebra}
(or {\it the trigonometric Cherednik algebra})
$\daffH$ of type $\GL_n$
is the associative $\F$-algebra which is generated
by the algebras $\F P$, $\F W$ and $S(\h)$, and subjects
to the following defining relations:
\begin{align*}
&s_i h =s_i (h) s_i-\bra\al_i\mid h\ket
\ (i\in[1,n-1],\ h\in\h),\quad
s_i e^\eta s_i= e^{s_i (\eta)}\
(i\in[1,n-1],\  \eta\in P),\\ 
&\left[ h, e^\eta\right] =\kappa\, \bra \eta\mid h\ket e^\eta
+\sum_{\al\in R^+}\bra \al\mid h\ket
{\frac{e^\eta-e^{s_\al(\eta)}}{ 1-e^{-\al}} s_{\al}}\ \
\left(h\in\h,\,\eta\in P\right),
\end{align*}
where $e^\eta$ denotes the element of $\F P$
corresponding to $\eta\in P$.

It is known that
$\F P\* \F W \* S(\h)\cong \daffH$ via the multiplication map,
and
the subalgebra $\F\W \cdot S(\h)\subset \daffH$ is isomorphic to the
degenerate affine Hecke algebra $\affH$.
Let $\lm,\mu\in\Pg^+$ such that $\lm-\mu\models n$.
Let $\sp_\lsm$ denote the irreducible module of $\affH$
associated with the skew diagram $\lsm$ (\cite{Ra,Ch;special_bases}).
Set $\daffst(\lm,\mu) 
=\daffH\otimes_{\affH}\sp_{\lsm}$.

When 
$\kappa\in\Z_{\geq0}$ and 
$\lm,\mu\in \Pg^+(\kappa)$, it is known that
the $\daffH$-module $\daffst(\lm,\mu)$ has
a unique simple quotient (\cite{AST}),
which we will denote by $\daffL(\lm,\mu)$.
%
\section{The degenerate double affine Hecke algebra and conformal field 
theory}
For an algebra $\alg$ and an $\alg$-module $M$, we
write $M_\alg=M/\alg M$.

Let $\kappa\in \F\setminus\{0\}$. 
Let $N$ be a highest weight $\affg$-module of level $\kappa-\m$
and let $M$ be a lowest weight $\affg$-module of level $-\kappa+\m$.
Introduce the following space of coinvariants:
\begin{align*}
  \daffC(N,M)
&=\left(N\otimes
\Vbox\xpm 1\otimes\dots\otimes \Vbox\xpm{n}\otimes M
\right)_{\g[t,t^{-1}]}.
\end{align*}
%
\begin{remark}
Write $\F[\underline{x}^{\pm1}]$ for the ring 
$\F[x_1^{\pm1},\dots,x_n^{\pm1}]$
of Laurent polynomials, which
can be identified with  the coordinate ring of the affine variety 
$\T=(\F\setminus\{0\})^n$.
Let
$\Xreg$ denote the localization
of
$\F[\underline{x}^{\pm1}]=\F[\T]$
at the diagonal set
$\triangle=\cup_{i<j}
\{(x_1,\dots,x_n)\in\T\mid x_i=x_j\}$,
that is,
$\Xreg= \F\left[x_1^{\pm1},\dots,x_n^{\pm1},\frac{1}{x_i-x_j}\
(i\leq j)\right].$
The space  $\daffC(N,M)$ can be regarded 
as an $\F[\underline{x}^{\pm1}]$-module. 
It follows that
the $\Xreg$-module
$\Xreg\*_{\F [\underline{x}^{\pm1}]}\daffC (N,M)
$ admits
an integrable connection,
and
it gives a vector bundle
on the affine variety $\T\setminus\triangle$ (\cite{AST,VV}).
In particular, for 
$\kappa\in\Z_{\geq0}$ and 
$\lm,\mu\in \Pg^+(\kappa)$, it can be shown that
the bundle $\Xreg\*_{\F [\underline{x}^{\pm1}]}\daffC
(\affLl(\mu),\affLl^\dagger(\lm))$
is equivalent to
the vector bundle of ``conformal coinvariants''
(or conformal blocks)
whose fiber at $(\xi_1,\dots,\xi_n)\in\T\setminus\triangle$
is given by
$
( \affLl(\mu)\*
\affLl(\fwt_1)^{\*n}\*\affLl(\lm^\dagger))_{\g({\mathbb P}^1\setminus
\{0,\xi_1,\dots,\xi_n,\infty\})}.$
Here  $\fwt_1\in\hg^*$ is the highest weight of $\Vbox$,
$\lm^\dag=-w_0(\lm)$ with $w_0$ being the longest element of $\W$,
and $\g( {\mathbb P}^1\setminus\{0,\xi_1,\dots,\xi_n,\infty \} )$
denotes the Lie algebra of $\g$-valued algebraic functions
on $ {\mathbb P}^1\setminus\{0, \xi_1,\dots,\xi_n,\infty \}$,
which acts on
$ \affLl(\mu)\*
\affLl(\fwt_1)^{\*n}\*\affLl(\lm^\dagger)$ through the
Laurent expansion at each point.
(See e.g.  \cite{BK,FJKLM} for precise definitions.)
\end{remark}
In \cite{AST},
an action of the algebra $\daffH$ on the space $\daffC(N,M)$
was constructed through
the Knizhnik-Zamolodchikov connection,
and some conjectures were proposed
 concerning the case where
$N$ and $M$ are integrable, that is,
$N=\affLl(\mu)$,
$M=\affLl^\dagger(\lm)$ for 
$\lm,\mu\in\Pg^+(\kappa)$.
More precisely, in \cite{AST}, the $\daffH$-module 
$\daffC(\affLl(\mu),\affLl^\dagger(\lm))$ 
is conjectured to be isomorphic to
$\daffL(\lm,\mu)$ when $\lm,\mu\in\Pg^+(\kappa)$
and $\lm-\mu\models n$ (\cite[Conjecture 5.4.2]{AST}),
and moreover,
a conjectural formula concerning 
a decomposition of the space
$\unip\daffC(\affLl(\mu),\affLl^\dagger(\lm))$,
where $\unip=\frac{1}{n!}\sum_{w\in W}w$,
into the weight spaces with respect to
the algebra $\unip S(E)\unip$ was proposed (\cite[Conjecture 6.2.6]{AST}).

The irreducible $\daffH$-modules
$\daffL(\lm,\mu)$ $(\lm,\mu\in\Pg^+(\kappa))$
are studied in \cite{SV} and their structures are
described explicitly in terms of tableaux on periodic skew diagrams.
Through this combinatorial description, 
the conjectures are proved.
\begin{theorem}\label{th;daff}
Conjecture 5.4.2 and Conjecture 6.2.6 in \cite{AST} are true.
\end{theorem}
The details of the proof will be published in the forthcoming papers.
In the rest of this note,
we will present the analogous statements
for the rational Cherednik algebra $\ratH$.
These statements are derived from
Theorem~\ref{th;daff}
through the results established in \cite{S;ratandtrig}, where
it is shown that $\ratH$ can be embedded into $\daffH$
as a subalgebra (see Remark~\ref{rem;embedding}),
and the induction functor
$\daffH\otimes_{\ratH}(-)$ gives a nice correspondence
between the category of $\ratH$-modules and
that of $\daffH$-modules.
%
\section{Space of coinvariants and the rational Cherednik algebra}
%
Let $\kappa\in \F\setminus\{0\}$.
For a lowest weight $\affg$-module $M$ of level $-\kappa+\m$,
we set
\begin{align*}
\tilratC(M)&=
\Vbox[x_1] \otimes\dots\otimes \Vbox[x_n]\* M,\\
 \ratC(M)&=\tilratC(M)_{\g[t]}
=\tilratC(M)/\g[t] \tilratC(M).
\end{align*}
Let $\sigma_{ij}\in \End_\F
(\F[x_1,\dots,x_n])$
denote the permutation of $x_i$ and $x_j$.
Let $\pi_{ij}\in\End_\F(\Vbox^{\*n})$ denote the permutation of
$i$-th and $j$-th component of the tensor product.
Note that
$\tilratC(M)\cong
\F[x_1,\dots,x_n]\* \Vbox^{\*n}\*  M$
as a space, through which
we regard 
$\sigma_{ij}$ and $\pi_{ij}$ as
elements in $\End_\F(\tilratC(M))$, and
we let $\F[x_1,\dots,x_n]$ act on $\tilratC(M)$
via multiplication.

Put $e_{kl}[p]=e_{kl}\* t^p\in\affg$ $(k,l\in[1,\m],\ p\in\Z)$,
where $e_{kl}\in\g$ denotes the matrix unit as before.
For $i\in[1,n]$, put
\begin{equation*}
\omega_i=\sum_{p\in\Z_{\geq 1}}\sum_{k,l\in[1,\m]}
1^{\*i-1}\* e_{kl}[p-1]\*1^{\*n-i}\* e_{lk}[-p], 
\end{equation*}
which is an element of
some completion of $U(\affg)^{\*n+1}$
and defines a well-defined operator on $\tilratC(M)$.
Define the operators on $\tilratC(M)$ by
\begin{equation*}
D_i=\kappa\frac{\partial}{\partial x_i}
+\sum_{j\in[1,n],\,j\neq i}
\frac{1}{x_i-x_j}
(1-\sigma_{ij})\pi_{ij}
-\omega_i
\quad (i\in[1,n]). 
\end{equation*}
A parallel argument as in  \cite[\S4]{AST} implies the following:
\begin{theorem}
Let $M$ be a lowest weight $\affg$-module of level $-\kappa+m$.

\smallskip\noindent
{\rm (i)} There exists a unique algebra homomorphism
$\theta:\ratH\to \End_\F(\tilratC(M))$
such that
\begin{align*}
\theta(\ech_i)&=D_i\  (i\in[1,n]),\quad
\theta(\e_i)=x_i\  (i\in[1,n]),\quad
\theta(s_i) =\pi_{i\, i+1}\sigma_{i\, i+1}\  (i\in[1,n-1]).
\end{align*}
\noindent
{\rm (ii)} The $\ratH$-action on $\tilratC(M)$ above preserves
the subspace $\g[t]\tilratC(M):$
$$\theta(a)\g[t]\tilratC(M)\subseteq\g[t]\tilratC(M)$$ for
all $a\in \ratH$.
Therefore, $\theta$ induces
an $\ratH$-module structure on $\ratC(M)$.
\end{theorem}
Define the elements $\u_i\ (i\in[1,n])$ of $\ratH$
by
$\u_i=\e_i\,\ech_i+\sum\nolimits_{j\in[1,i-1]}s_{ji}.$
It follows that $\u_1,\dots,\u_n$ are pairwise commutative.
Write $\yy$ for the subalgebra $\F[\u_1,\dots,\u_n]$ of
$\ratH$.
\begin{remark}
\label{rem;embedding}
As pointed out in \cite{S;ratandtrig},
there exists an algebra embedding $\emb:\ratH\to\daffH$
such that
\begin{equation*}
\emb(s_i)= s_i\ (i\in[1,n-1]),\ \
\emb(\e_i)=e^{\e_i},\ \  
\emb(\ech_i)=
e^{-\e_i}\left(\ech_i-\sum_{j\in[1,i-1]}s_{ji}\right)\ (i\in[1,n]).
\end{equation*}
The elements $\u_1,\dots,\u_n$ of $\ratH$ correspond
to the elements $\ech_1,\dots,\ech_n$ of $\daffH$ via this embedding,
and the subalgebra $\F W\cdot \yy$ of $\ratH$
is isomorphic to the degenerate affine Hecke algebra $\affH$. %
\end{remark}

A finitely generated $\ratH$-module $M$ is said to be {\it calibrated}
if $M$ is locally nilpotent for $\ech_1,\dots,\ech_n$ and 
$M=\oplus_{\zeta\in \F^n}M_\zeta$,
where $M_\zeta=\{v\in M\mid u_iv=\zeta_i v\ \forall i\in[1,n]\}$
with $\zeta=(\zeta_1,\dots,\zeta_n)$.

Now,
let $\kappa\in\Z_{\geq1}$ 
and
set $\Lambda^+_{\kappa}(\m,n)=\Pg^+(\kappa)\cap \Lambda^+(\m,n)$, where
 $\Lambda^+(\m,n)=\{\lm\in\Pg^+\mid \lm\models n,\
\lm_\m\geq1\}$.
It is shown in \cite[Theorem 7.2]{S;ratandtrig}
that the assignment $\lm\mapsto \ratL(\lm)$
induces a bijection between the set $\sqcup_{\m\in[1,\kappa]}
\Lambda^+_{\kappa}(\m,n)$ and the
set of isomorphism classes of
the irreducible calibrated $\ratH$-modules. 

As a rational analogue of 
\cite[Conjecture~5.4.2]{AST}, we have
\begin{theorem}\label{th;image_rat}
Let $\kappa, \m\in \Z_{\geq1}$
and $\lm\in \Lambda_{\kappa}^+(\m,n)$.
Then
$\ratC(\affLl^\dagger(\lm))\cong \ratL(\lm)$
as an $\ratH$-module.
\end{theorem}
\begin{remark}
As a consequence of Theorem~\ref{th;image_rat} and 
\cite[Theorem 7.2]{S;ratandtrig},
any irreducible calibrated $\ratH$-module 
is realized as the space $\ratC(\affLl^\dagger(\lm))$ of coinvariants
for some $\lm\in\sqcup_{\m\in[1,\kappa]}\Lambda^+_{\kappa}(\m,n)$.
\end{remark}
%
\section{Character formula and Kostka polynomials} 
%
Put $\unip=\frac{1}{n!}\sum_{w\in W}w\in \F W$.
For an $\ratH$-module $M$, 
the space $\unip M$ is called the spherical part (or the symmetric part)
of $M$.
Recall that $\ratH$ has the canonical grading operator
 $\der
=\kappa^{-1}(\sum_{i=1}^{n}\e_i\ech_i+\sum_{i<j}s_{ij})
=\kappa^{-1}\sum_{i=1}^{n}\u_i\in\yy$,
which
satisfies $[\der,{\e_i}]=\e_i,\ [\der,\ech_i]=-\ech_i$ for $i\in[1,n]$
and $[\der,w]=0$ for $w\in\W$.
Observe that $\der$ preserves the spherical part of any $\ratH$-module.

Let $\lm\in\Lambda^+_\kappa(\m,n)$
and put
$\unip\ratL(\lm)_{(k)}=\{v\in \unip\ratL(\lm)\mid\der v=kv\}$ 
for $k\in\F$.
Then it follows
 that $\dim_\F \unip\ratL(\lm)_{(k)}<\infty$ and
 $\unip\ratL(\lm)=\oplus_{k\in\Z+\cas_\lm}\unip\ratL(\lm)_{(k)}$,
 where
$$\cas_\lm=\frac{1}{2\kappa}\sum_{i=1}^\m
\lm_i(\lm_i-2i+1).$$
Hence we have
$\Ch{\unip\ratL(\lm)}=\sum_{k\in\Z+\cas_\lm}q^k
\dim_\F \unip\ratL(\lm)_{(k)}$.

Analogously to \cite[Corollary 6.2.6]{AST},
a decomposition of $\unip\ratL(\lm)$ into the weight spaces
with respect to the algebra $\unip\yy \unip$ is described explicitly.
This gives
a simple and remarkable formula for $\Ch{{\unip\ratL(\lm)}}$,
 which we will describe below.

Let $\lm\in\Lambda^+(\m,n)$.
Let $T$ be a standard tableau on the Young diagram
$\lm$; namely, $T$ is a bijection from the diagram
$\lm$ (viewed as a subset of $\Z^2$, see \eqref{eq;youngdiagram})
to the set $[1,n]$
satisfying $T(a+1,b)>T(a,b)$ when $(a,b),(a+1,b)\in\lm$,
and
 $T(a,b+1)>T(a,b)$ when $(a,b),(a,b+1)\in\lm$.
Define the sequence 
$\emptyset=\lm_T^{(0)}\subset \lm_T^{(1)}\subset \dots
\subset  \lm_T^{(n)}=\lm$
of subdiagrams of $\lm$ by
$\lm_T^{(i)}=T^{-1}[1,i]$.
For $\l\in\Z_{\geq0}$,
a tableau  $T$ is called an
$\l$-{\it restricted standard tableau}
if $\lm_T^{(i)}\in \Pg^+(\l+\m)$ for all $i\in [1,n]$.
Let
${\mathrm St}^{(\l)}(\lm)$ denote the set of 
$\l$-{restricted} standard tableaux
on $\lm$.
Observe that ${\mathrm St}^{(\l)}(\lm)=\emptyset$
unless $\lm\in\Lambda^+_{\l+\m}(\m,n)$.
For a standard tableau $T$ 
 and $i\in[1,n]$,
define
 \begin{equation}
 h_i(T)=\begin{cases}
 &1\quad \hbox{if}\ a<a',\\
 &0\quad \hbox{if}\ a\geq a',
 \end{cases}
\ \hbox{ where }(a,b)=T^{-1}(i),\ (a',b')=T^{-1}(i+1),
 \end{equation}
 and put $\charge(T)={\sum_{i=1}^{n}(n-i)h_i(T)}.$
The map $d$ is called the {\it cocharge}, and it is 
 essentially the same with
what is called the energy function in the RSOS model
(see e.g. \cite{NY}).
Put
\begin{equation} \label{eq;F}
{\check K}_{\lm\; (1^n)}^{(\l)}(q)=
\sum_{T\in\St^{(\l)}(\lm)}q^{\charge(T)},
\end{equation}
where $(1^n)=(1,1,\dots,1)\in \Lambda^+(n,n)$.
Note that the polynomial given by
$
K_{\lm\; (1^n)}^{(\l)}(q)=
q^{\frac{1}{2} n(n-1)}
 {\check K}_{\lm\; (1^n)}^{(\l)}(q^{-1})
$
is called the $\l$-restricted Kostka polynomial
associated with the partitions $\lm$ and $(1^n)$
(see e.g. \cite{SS,FJKLM}).
%
\begin{theorem}\label{th;char_rat}
Let $\m,\kappa \in\Z_{\geq1}$
and
 $\lm\in \Lambda^+_{\kappa}(\m,n)$ .
Then
  \begin{equation*}
    \Ch{\unip\ratL(\lm)}=
\frac{q^{\cas_\lm}}{\mathop{\Pi}\limits_{i=1}^n(1-q^i)}
{\check K}_{\lm\; (1^n)}^{(\kappa-\m)}(q)
=
\frac{q^{\cas_\lm+\frac{1}{2}n(n-1)}}{\mathop{\Pi}\limits_{i=1}^n
(1-q^i)}
{K}_{\lm\; (1^n)}^{(\kappa-\m)}(q^{-1}).
  \end{equation*}
\end{theorem}
\begin{remark}
In \cite[Theorem 1.10]{BEG}, the character 
$\Ch{\unip\ratL(\lm)}$ was calculated when
$\lm=(n)$ and $\kappa=-\frac{n}{r}$,
 where $r\in\Z_{\geq1}$ with $(n,r)=1$.
\end{remark}

\end{document}